\documentclass[leqno,11pt,twoside]{article}
\usepackage{mathrsfs}
\usepackage{amssymb}
\long\def\symbolfootnote[#1]#2{\begingroup%
\def\thefootnote{\fnsymbol{footnote}}\footnote[#1]{#2}\endgroup}

\makeatletter
   \renewcommand{\section}{\@startsection{section}{1}{0mm}
   {\baselineskip}%
   {\baselineskip}{\normalfont\normalsize\bf
	\centering}}%
   \makeatother

\def\thebibliography#1{\footnotesize\par\vskip4mm plus 1mm minus 1mm
\centerline{\normalsize\textbf{References}}
\par\nobreak\list
 {[\arabic{enumi}]}{\settowidth\labelwidth{[#1]}
 \leftmargin\labelwidth\advance\leftmargin\labelsep
 \partopsep 0pt \itemsep 0pt \parsep 0pt
 \usecounter{enumi}}
 \def\newblock{\hskip .11em plus .33em minus -.07em}
 \sloppy
 \sfcode`\.=1000\relax}

\newtheorem{definition}{Definition}

\newtheorem{prop}[definition]{Proposition}
\newtheorem{thm}[definition]{Theorem}
\newtheorem{cor}[definition]{Corollary}
\newtheorem{rmk}[definition]{Remark}

\def\CC{{\mathscr C}}

\def\FF{{\mathscr F}}
\def\GG{{\mathscr G}}

\def\JJ{{\mathscr J}}

\def\MM{{\mathscr M}}

\def\OO{{\mathscr O}}

\def\UU{{\mathscr U}}

\def\XX{{\mathscr X}}
\def\YY{{\mathscr Y}}
\def\ZZ{{\mathscr Z}}

\def\va{\mbox{\rm va}}

\def\Pic{\mathop{\rm Pic}\nolimits}

\def\deg{\mathop{\rm deg}\nolimits}
\def\ker{\mathop{\rm ker}\nolimits}

\def\dim{\mathop{\rm dim\, }\nolimits}

\def\Cliff{\mathop{\rm Cliff}\nolimits}

\def\proof{
  \noindent
  {\bf Proof.}
}

\def\proofmain1{
  \noindent
  {\bf Proof of Theorem \ref{thm: main1}.}
}

\def\proofodd{
  \noindent
  {\bf Proof of Theorem \ref{thm: main2}.}
}

\def\endproof{
{\unskip\nobreak\hfil\penalty50\hskip2em\hbox{}\nobreak\hfill
          $\square$\bigbreak}
}

\begin{document}

\centerline{\bf REMARKS ON SYZYGIES OF $d$-GONAL CURVES}

\bigskip

\centerline{\sc Marian Aprodu}

\vskip 1cm

\begin{minipage}{12.5cm}\small {\sc Abstract.}
We apply a degenerate version of a result due to Hirschowitz, Ramanan 
and Voisin to verify Green and Green-Lazarsfeld conjectures over explicit 
open sets inside each $d$-gonal stratum of curves $X$ with $d<[g_X/2]+2$. 
By the same method, we verify the Green-Lazarsfeld conjecture for 
any curve of odd genus and maximal gonality. The proof invokes 
Voisin's solution to the generic Green conjecture as a key argument.
\end{minipage}

\symbolfootnote[0]{2000 {\em Mathematical Subject Classification}.
14H51, 14C20, 13D02.}
\symbolfootnote[0]{{\em Key words and phrases}: syzygies, curve, 
gonality, Clifford index, stable curve.}
\symbolfootnote[0]{This work was supported by a Humboldt 
Research Fellowship.}

\section{Introduction}

A result due to Green and Lazarsfeld, {\em cf.} 
\cite[Appendix]{G}, shows that the existence of special linear systems on 
a complex projective variety $X$ reflects into non-vanishing of certain 
Koszul cohomology groups $K_{p,1}$. Recall that for two integers
$p$ and $q$ and a line bundle $L$ on $X$, the Koszul cohomology group 
$K_{p,q}(X,L)$ was defined in \cite{G} as the cohomology of the complex:
$$
\wedge^{p+1} H^0(L)\otimes H^0(L^{q-1})
\rightarrow
\wedge^p H^0(L)\otimes H^0(L^q)
\rightarrow
\wedge^{p-1} H^0(L)\otimes H^0(L^{q+1}).
$$

On a smooth, connected curve $X$, one expects that the Green-Lazarsfeld 
non-vanishing result be revertible, that is non-trivial $K_{p,1}$ groups
with values in line bundles chosen in convenient ways would give 
rise to special linear systems on the curves. This amounts to proving 
that particular Koszul cohomology groups vanish if some
linear series loci are empty.

For the canonical bundle, Green conjectured that $K_{g-c-1,1}(X,K_X)=0$ 
for any curve $X$ of genus $g$ and Clifford index $c$, {\em see} 
\cite{G}. Similarly, if $X$ is of gonality $d$, Green and Lazarsfeld 
predicted $K_{h^0(L)-d,1}(X,L)=0$ for any line bundle $L$ of
sufficiently large degree, {\em see} \cite{GL1}. Put in other words,
the two basic invariants of a curve, Clifford index and gonality,
would be read off suitable Koszul cohomology groups. It is important
to mention that they are related by the inequalities 
$d-3\le c\le d-2\le [(g_X-1)/2]$, 
{\em cf.} \cite{ACGH} and \cite{CM},

Both Green, and Green-Lazarsfeld conjectures have been verified for 
general curves of any genus (for Green's conjecture {\em see} \cite{V2} 
and \cite{V3}, and for the Green-Lazarsfeld conjecture we refer to 
\cite{AV} and \cite{A2}) with the notable difference that in the odd-genus  
case, the Green conjecture is known to hold for any curve of 
maximal gonality, as shown by Hirschowitz, Ramanan and 
Voisin, {\em cf}. \cite{HR} and \cite{V3}. 

\begin{thm}[Hirschowitz-Ramanan-Voisin]
\label{thm: HRV}
Any smooth curve $X$ of genus $2k+1$ with $K_{k,1}(X,\omega_X)\ne 0$
carries a pencil of degree $k+1$.
\end{thm}

The proof of this remarkable result relies 
on a comparison of two divisors in the moduli space $\MM_{2k+1}$, the one 
of curves which carry a $g^1_{k+1}$ and, respectively, the locus of 
canonical curves with extra-syzygies. Already showing that the latter one 
is a genuine divisor, and does not cover the whole moduli space is 
hard work, and was concluded only recently in Voisin's tour de force 
\cite{V2}, \cite{V3}, after having 
developed a completely new approach to the problem in \cite{V2}.

\medskip

The Green conjecture is known to hold also for curves of non-maximal
gonality which are generic in their gonality strata, and this
happens for all possible gonalities {\em cf.} \cite[Corollary 1]{V2} 
and \cite[Theorem 0.4]{T}, {\em see} also \cite[Theorem]{S1}. 
It is well-known that by fixing the gonality $d$ we obtain a stratification 
of the moduli space of curves with irreducible gonality strata, the maximal 
possible value for $d$ corresponding to the open stratum. Therefore, 
it makes perfect sense to speak about general points in a  given gonality 
stratum. 

Starting from the generic Green conjecture, Voisin had the idea to 
degenerate smooth curves on a $K3$ surface to irreducible nodal curves,
in the same linear system, in order to make the Koszul cohomology of the 
latter ones vanish, and so she verified the Green conjecture for the 
normalizations of these curves, {\em cf.} \cite{V2}. This fact resulted 
into a very short and elegant solution for the Green conjecture for 
generic curves $X$ of non-maximal gonality larger than $g_X/3$. Then she 
observed furthermore that exactly the same method yields, in completion 
to \cite[Theorem 1]{A1}, to a solution for the Green-Lazarsfeld conjecture 
for generic curves of fixed gonality in the same order range, {\em cf.} 
\cite[Theorem 1.3]{AV} and \cite[Theorem 1.4]{AV}. This strategy did much 
better than the partial \cite[Theorem 3]{A1}, {\em see} \cite[Remark 
2]{AV} for some related comments. 

At the other end of the spectrum, normalizations of irreducible nodal 
curves on $\mathbb{P}^1\times\mathbb{P}^1$ have been used to verify the 
two conjectures for generic curves $X$ in any stratum of gonality
bounded in the order $\sqrt{g_X}$, {\em cf.} \cite[Theorem]{S1}, and 
\cite[Theorem 4]{A1}, the vanishing of the Koszul cohomology having been 
proved, in contrast to Voisin, directly on the normalizations.
Further connections between Green and Green-Lazarsfeld conjectures
are emphasized in \cite[Appendix and II]{A1}.

Having had these indications pleading for unity, it seems natural to try 
to find a common space of curves with fixed gonality on which the two 
conjectures be treated in a unitary way -- this is the goal of the present 
work. We use degenerations, taking the path shown by Voisin, 
which is to compute the Koszul cohomology directly on nodal curves, instead 
of normalizing them first. The proof of our first result will show that 
this is a very natural thing to do.

\begin{thm}
\label{thm: main1} 
Let $d\ge 3$  be an integer, and $X$ be a smooth $d$-gonal curve 
with $d<[g_X/2]+2$, and such that $\dim(W^1_{d+n}(X))\le n$ for all 
$0\le n\le g_X-2d+2$. Then $\Cliff(X)=d-2$, and $X$
verifies both Green, and Green-Lazarsfeld 
conjectures. 
\end{thm}

In the statement, $W^1_{d+n}(X)$ denotes the subvariety of 
$\Pic_{d+n}(X)$ of bundles with two or more independent sections.
We mention that, since $X$ is non-hyperelliptic, we always have the bound
$\dim(W^1_{d+n}(X))\le d+n-3$ for all $0\le n\le g_X-2d+2$, {\em see}
\cite{HMa}. Beside, according to a problem raised by G. Martens, {\em cf.}
\cite[Statement (T), p. 280]{GMa}, it is very likely that the conditions 
appearing in the statement of Theorem \ref{thm: main1}
be reduced to the single condition $\dim(W^1_{d}(X))=0$.
The locus of curves satisfying the hypotheses of the theorem is anyway
non-empty and open in the $d$-gonal stratum, and all but one possible 
gonality strata are concerned.
In particular, Theorem \ref{thm: main1} gives an alternate proof to 
\cite[Theorem 0.4]{T}, \cite[Theorem]{S1}, \cite[Theorem 4]{A1}, and,
for the Green-Lazarsfeld conjecture, it fills in the existing gap 
for generic $d$-gonal curves with $d$ in 
order between $\sqrt{g_X}$ and $g_X/3$.

For small $d$ one can specify the curves which do not satisfy the 
hypothesis, {\em cf.} \cite{ACGH}, \cite{HMa}, \cite{Mu}, \cite{Ke}. 
Among the new effective results that follow directly from Theorem 
\ref{thm: main1} we count then the validity of the Green-Lazarsfeld 
conjecture for non-bielliptic tetragonal curves, and of the two 
conjectures for pentagonal and for hexagonal curves which are neither
trielliptic nor double coverings of curves of genus three, 
{\em see} Corollary \ref{cor: left out} in Section 3.

The extremal case $g_X-2d+2=0$ indicates that both Green, and 
Green-Lazarsfeld conjectures are verified for any curve $X$ of
even genus and maximal gonality which carries finitely many
minimal pencils. From the Green conjecture point of view, this 
can be seen as a first step in finding a correspondent of the 
Hirschowitz-Ramanan-Voisin Theorem in the even-genus case. 
Other consequence of Theorem \ref{thm: main1} regards curves of 
even genus and maximal Clifford index on $K3$ surfaces, {\em cf.} 
Corollary \ref{cor: k3}.

\medskip

The case not covered by Theorem \ref{thm: main1} is, as expected, the
case of curves of odd genus and maximal gonality. The second result, 
which refines \cite[Theorem 1.1]{A2}, is concerned with the 
Green-Lazarsfeld conjecture for these curves, for which the Green 
conjecture has already been settled.

\begin{thm}
\label{thm: main2} The Green-Lazarsfeld conjecture is valid for any 
smooth curve $X$ of genus $g_X=2k-1$ and gonality $k+1$, with $k\ge 2$.
\end{thm}

We obtain directly from Theorem \ref{thm: main2}
a version of the Hirschowitz-Ramanan-Voisin 
Theorem for syzygies of pluricanonical curves.

\begin{cor}
\label{cor: pluricanonical}
If $X$ is a  smooth curve of genus $2k-1$, where $k\ge 2$, with
$$K_{2(2p-1)(k-1)-(k+1),1}(X,\omega_X^{\otimes p})
\ne 0$$
for some $p\ge 2$, then $X$  carries a $g^1_k$.
\end{cor}

Roughly speaking, the proof strategy of the two theorems above is the same, 
namely to construct, starting from $X$, some suitable singular stable 
curves of genus $g=2k+1$ and analyze their syzygies. A degenerate version 
of the Hirschowitz-Ramanan-Voisin Theorem, proved in Section 2, 
tells us that if these stable curves have extra-syzygies, then 
they are limits of smooth curves with extra-syzygies, and this fact
translates into the existence of special torsion-free sheaves 
of rank one, which will trace pencils on the original
curve $X$. For Theorem \ref{thm: main1}, the singular curves which 
we construct are irreducible, and this choice was inspired directly 
by the ideas of Voisin applied in \cite{V2} and \cite{AV}, and which 
we have briefly explained above. This result is then a natural 
continuation of \cite[Corollary 1]{V2} and \cite[Theorems 1.3
and 1.4]{AV}. 
The difference is that we work now with abstract families of 
curves degenerating to a stable curve, whereas in \cite{V2} and 
\cite{AV}, these degenerations were always made on $K3$ surfaces. 
The freedom we gain also
reflects into the possibility of working with curves with arbitrarily many 
nodes, getting rid of the main obstruction to improving
the statements quoted above. As for Theorem \ref{thm: main2},
the singular curves which we construct are curves with one smooth 
rational component, similarly to \cite[Theorem 1.1]{A2}. 
These curves are themselves limits of irreducible curves 
with two nodes, indicating this case as being {\em "more degenerate"}
than the others, and giving a hint of why should it be treated
separately. The origin of this difference is that this
is the only case of a gonality stratum whose general members
carry infinitely many minimal pencils.

\section{Syzygies of singular stable curves with very ample 
canonical bundle}

Throughout this section the arithmetic genus will remain fixed to 
$g=2k+1$, where $k\ge 2$. We study the Koszul cohomology of a
stable curve $Y$ of genus $g$ in relation with its geometry.
Recall that {\em stable curves} are reduced connected curves with finite 
group of automorphisms, and with only simple double points (nodes) as 
possible singularities. They have been introduced by Deligne and Mumford 
with the aim of compactifying the moduli space $\MM_g$ of smooth curves of 
genus $g$. Singular stable curves of arithmetic genus $g$ lie on a 
normal-crossing divisor $\Delta_0\cup\dots\cup\Delta_{[g/2]}$ in 
$\overline{\MM_g}$, on the boundary of $\MM_g$, 
and the general element in $\Delta_0$ is 
irreducible, whereas a general element in $\Delta_i$
is the union of two curves of genus $i$ and $g-i$ respectively,
meeting in one point. 

\medskip

From now on, we shall work exclusively with stable curves with 
{\em very ample canonical bundle}, {\em see} \cite[Theorem F]{C} 
and \cite[Theorem 3.6]{CFHR} for precise criteria. 

\medskip

\noindent
{\bf Notation.}
The open subspace in $\overline{\MM_g}$ of points corresponding 
to stable curves of genus $g$ with very ample canonical bundle, 
which can be easily shown to be contained in $\MM_g\cup\Delta_0$, 
will be denoted by $\MM^{\va}_g$. 

\medskip

Let $[Y]\in \MM^{\va}_g$, 
and denote, for simplicity $\mathbb{P}:=\mathbb{P}H^0(Y,\omega_Y)^*$, 
which contains the image of $Y$, set $Q=T_{\mathbb{P}}(-1)$ the universal 
quotient bundle, and $Q_Y$ the restriction of $Q$ to $Y$. The Koszul 
cohomology of $Y$ with values in $\omega_Y$ has then the following 
description, \cite[Proposition 2.1]{HR}.

\begin{prop}[Hirschowitz-Ramanan]
With the notation above, for any $n\ge 1$, there exists an isomorphism
$$
K_{n,1}(Y,\omega_Y)\cong\ker\left(H^0(\mathbb{P},\wedge^{2k-n+1}Q(1))
\rightarrow H^0(Y,\wedge^{2k-n+1}Q_Y\otimes\omega_Y)\right).
$$
\end{prop}

\begin{rmk}
\label{rmk: VB}
{\rm 
For the choice $n=k$ and for a smooth curve $Y$, Hirschowitz and Ramanan 
remarked that the two spaces appearing in the description above, namely
$H^0(\mathbb{P},\wedge^{k+1}Q(1))$ and 
$H^0(Y,\wedge^{k+1}Q_Y\otimes\omega_Y)$, have the same dimension,
{\em cf.} \cite[Proof of Proposition 3.1]{HR}. 
This observation was essential in the proof of their main result.
Applying the Riemann-Roch Theorem, and the Serre duality we observe 
that the two spaces in question are still of the same dimension even if
$Y$ was a singular stable curve with very ample canonical bundle, 
since $H^0(Y,\wedge^{k+1}Q_Y^*)=0$.
The vanishing of $H^0(Y,\wedge^{m+1}Q_Y^*)$ for any $m\ge 0$ is a standard 
fact, and is implied by the following classical remarks. Firstly, we 
know that there are natural exact sequences, for any $p$ and $q$,
$$
0\rightarrow H^0(Y,\wedge^{p+1}Q_Y^*\otimes\omega_Y^{q-1})
\rightarrow 
\wedge^{p+1}H^0(Y,\omega_Y)\otimes H^0(Y,\omega_Y^{q-1})
\rightarrow H^0(Y,\wedge^pQ_Y^*\otimes\omega_Y^q).
$$
Secondly, the Koszul differential $\wedge^{m+1}H^0(\omega_Y)\rightarrow 
\wedge^mH^0(\omega_Y)\otimes H^0(\omega_Y)$ is injective, and it factors 
through the inclusion of $H^0(\wedge^mQ_Y^*\otimes\omega_Y)$ into 
$\wedge^mH^0(Y,\omega_Y)\otimes H^0(Y,\omega_Y)$.}
\end{rmk}

\noindent
{\bf Convention.} Hirschowitz and Ramanan used the term 
{\em with extra-syzygies} is used to designate a smooth curve $Y$ 
of genus $2k+1$ with $K_{k,1}(Y,\omega_Y)\ne 0$. We adopt this 
terminology and extend it to singular curves with the same 
non-vanishing property.

\medskip

From the syzygy point of view, singular stable curves with very ample
canonical bundle behave similarly to smooth
curves. For instance, those which have extra-syzygies, are 
degenerations of smooth curves with extra-syzygies.

\begin{prop}
\label{lemma: HRV}
Let $Y$ be a singular stable curve of genus $g=2k+1$ with very ample
canonical bundle. 
Then $Y$ has extra-syzygies if and only if $[Y]$ belongs to 
the closure of the locus of $(k+1)$-gonal smooth curves.
\end{prop}

\proof We set some notation first. Following \cite{HM}, let $D_{k+1}$ be 
the divisor on $\MM_g$ of curves with a pencil of degree $k+1$, let 
$\overline{D}_{k+1}$ be its closure in $\overline{\MM_g}$, 
$D^{\va}_{k+1}=\overline{D}_{k+1}\cap \MM_g^{\va}$, and
$\Delta_0^{\va}=\Delta_0\cap \MM^{\va}_g$. 

By semi-continuity, the locus of curves in $\MM^{\va}_g$ with 
extra-syzygies is closed. Similarly to \cite{HR}, we show that 
it is actually a divisor. It amounts to proving that its inverse 
image on a covering $S^{\va}\rightarrow\MM^{\va}_g$ is a 
divisor. We choose a smooth $S^{\va}$ on which an 
universal curve exists. From \cite[Proof of Proposition 3.1]{HR} and 
from Remark \ref{rmk: VB} it follows 
that the locus of points on $S^{\va}$ corresponding to 
curves with extra-syzygies is the degeneracy locus of a morphism of
vector bundles of the same rank. Then this locus is a divisor, since we 
know it is not the whole space, {\em cf}. \cite[Theorem 4]{V3}.

Back to $\MM^{\va}_g$, the computation of \cite{HR} shows that the 
divisor of curves with extra-syzygies is equal to a multiple of 
$D^{\va}_{k+1}$ plus, eventually, a multiple of $\Delta_0^{\va}$
(recall that $\Delta_0$ is irreducible, and so is $\Delta_0^{\va}$). 
The possibility that the whole $\Delta_0^{\va}$ be contained in the locus 
of curves with extra-syzygies is ruled out by the same result of \cite{V3}.
The example of a curve in $\Delta_0^{\va}$ with no extra-syzygies given by 
Voisin, and used in \cite[Proof of Theorem 1.4]{AV}, is an irreducible 
curve with one node lying on a $K3$ surface, {\em compare} to the proof 
of Theorem \ref{thm: main1} in the next section. 
In particular, it follows that a curve in $\Delta_0^{\va}$ 
has extra syzygies if and only if it belongs to $D_{k+1}^{\va}$, 
which we wanted to prove.
\endproof

We analyze next further consequences of having extra-syzygies, and show 
that this condition yields to the existence of certain suitable sheaves. 
By means of Proposition \ref{lemma: HRV}, a singular
stable curve $Y$ with $\omega_Y$ very ample and with extra-syzygies lies 
in an one-dimensional flat family $\CC\rightarrow T$ of curves such that 
$\CC_{t_0}\cong Y$, and $\CC_t$ are smooth, and belong to $D_{k+1}$ for 
$t\ne t_0$. For safety, let us make the further assumption that 
$Y$ has trivial automorphism group, altough this extra-condition
might be unnecessary. By shrinking $T$ if
needed, we can make the same hypothesis for the curves 
$\CC_t$. By the compactification theory of the generalized 
relative Jacobian, {\em see} \cite{Ca} and \cite{P}, there exists 
a family $\JJ_{1-k}(\CC/T)$, flat and proper over $T$, whose
fiber over $t\ne t_0$ is just the Jacobian variety of line 
bundles of degree $k+1$, whereas the fiber over $t_0$ parametrizes 
$gr$-equivalence classes of torsion-free, $\omega_Y$-semistable 
sheaves $F$ of rank one on $Y$ with $\chi(F)=1-k$; by definition,
a coherent sheaf on $Y$ is {\em torsion-free} if
it has no non-zero subsheaf with zero-dimensional support.
It follows that the subspace of pairs
$\{ (\FF_t,\CC_t)\in \JJ_{1-k}(\CC/T)\times_T\CC
,\; h^0(\CC_t,\FF_t)\ge 2\}$ is closed in the fibered product,
and, since $[\CC_t]\in D_{k+1}$ for all $t\ne t_0$, we are lead
to the following degenerate version of the Hirschowitz-Ramanan-Voisin
Theorem.

\begin{prop}
\label{prop: F}
Let $Y$ be a singular stable curve of genus $g=2k+1$ with very ample 
canonical bundle and trivial automorphism group. If $Y$ has
extra-syzygies, then there exists a torsion-free, $\omega_Y$-semistable
sheaf $F$ of rank one on $Y$ with $\chi(Y,F)=1-k$ and $h^0(Y,F)\ge 2$. 
\end{prop}

In several particular cases, this sheaf can be analyzed even further
to deduce the existence of other suitable linear systems
on the normalizations of components of $Y$, {\em see}
the proofs of Theorems \ref{thm: main1} and 
\ref{thm: main2}, and {\em compare to} \cite[Corollary 1, p. 68]{HM}.

\section{Proof of Theorem \ref{thm: main1} and consequences}

We recall the following result from \cite{A1}, which will be used for 
the Green-Lazarsfeld conjecture, {\em see}
also \cite[Theorem 2.1]{AV}:
\begin{thm}
\label{thm: Ap}
If $L$ is a nonspecial line bundle on a smooth curve $X$, which satisfies 
\linebreak $K_{n,1}(X,L)=0$, for an integer $n\ge 1$, then, for any 
effective divisor $E$ of degree $e\ge 1$, one has $K_{n+e,1}(X,L+E)=0$.
\end{thm}

\proofmain1
We set some notation. Define $k=g_X-d+1\ge 1$ and $\nu=g_X-2d+2\ge 0$, 
so that $d=k-\nu+1$, $g_X=2k-\nu$, and $0\le n\le\nu$.

\medskip

One should ensure firstly that curves satisfying the hypothesis of the
theorem exist. In the case $\nu=0$, this is automatic, but in the other 
cases it is less clear. We assume therefore that $\nu\ge 1$, condition
which is equivalent to $d<g_X/2+1$, and adopt the notation of \cite{AC}. 
From \cite[Theorem (2.6)]{AC}, we know
that a general $d$-gonal curve $X$ has precisely one minimal pencil, 
and moreover, any other base-point-free $g^1_h$ on $X$ with $h<g_X/2+1$ 
is composed with the given $g^1_d$. We count now the possible dimensions
of the irreducible components of the varieties $G^1_{d+n}(X)$ 
parametrizing the $g^1_{d+n}$'s on a generic $d$-gonal
curve $X$, and prove that they cannot exceed $n$. Since 
$G^1_{d+n}(X)$ is the canonical blowup of the $W^1_{d+n}(X)$ endowed 
with the determinantal structure, \cite[Chapter IV.3]{ACGH}, it would
prove our existence claim. 

Let $S$ be a smooth variety covering $\MM_{2k-\nu}$ on which there
exists an universal curve $\CC\rightarrow S$, let $\GG^1_d$, and 
$\GG^1_{d+n}$ the varieties over $S$ parameterizing pairs of smooth 
curves of genus $2k-\nu$, and pencils of degree $d$, and respectively 
$d+n$, and let $\XX$ be an 
irreducible component of $\GG^1_d\times_S\GG^1_{d+n}$ which
dominates $\GG^1_d$. Since $\GG^1_d$ is birational to the locus of 
$d$-gonal curves, {\em cf.} \cite[Theorem (2.6)]{AC} its dimension 
equals $2(2k-\nu)+2d-5$. 
Therefore, we must prove that $\dim(\XX)\le 2(2k-\nu)+2d+n-5$, which 
would ensure that a general fibre $\XX_\xi$ of the 
projection $\XX\rightarrow\GG^1_d$ is of dimension at most $n$. 

Let us suppose first that the general element of $\XX$ corresponds to a 
smooth $d$-gonal curve $X$ with a base-point-free $g^1_{d+n}$ which is 
not composed with a rational involution. Since $\XX$ dominates $\GG^1_d$, 
and minimal pencils on general $d$-gonal curves are simple, 
the hypotheses of \cite[Proposition (2.4)]{AC} are fulfilled. 
Coupled with \cite[Theorem (2.6)]{AC}, it implies that either $n=0$, 
and the two $g^1_d$'s for a generic point in $\XX$ coincide, 
in which case the projection from $\XX$ to $\GG^1_d$ is birational, or 
$d+n\ge (2k-\nu)/2+1$ and then
$\dim(\XX)=(2k-\nu)+2d+2(d+n)-7\le 2(2k-\nu)+2d+n-5$.

The other components have the following description.
Suppose that there is a non-empty open
subspace $\UU$ of $\XX$ whose points correspond to smooth 
$d$-gonal curves together with a $g^1_{d+n}$ having precisely 
$\mu$ base-points, 
where $0\le \mu\le n$. Then there is a natural map
$\UU\rightarrow \GG^1_d \times_s\GG^1_{d+n-\mu}$,
obtained by erasing the base points of the corresponding $g^1_{d+n}$'s.
Since $\UU$ is irreducible, its image is contained in an
irreducible component $\YY$ of $\GG^1_d\times_S\GG^1_{d+n-\mu}$,
whose generic point corresponds to a $d$-gonal curve together
with a base-point-free $g^1_{d+n-\mu}$ (recall that base-point-freeness
is an open condition). The projections on $\GG^1_d$ commute 
with the corestriction $\UU\rightarrow \YY$, so that $\YY$ dominates
$\GG^1_d$, too. If one denotes 
$\CC^{(\mu)}=(\CC\times_S\dots\times_S\CC)/\mathfrak{S}_\mu\rightarrow S$, 
where the product was taken $\mu$ times, 
the relative symmetric product, then the
induced map $\UU\rightarrow\YY\times_S\CC^{(\mu)}$
defined by splitting the corresponding $g^1_{d+n}$'s into their
free and fixed parts respectively, is injective. In
particular, $\dim(\XX)\le \dim(\YY)+\mu$. We prove next that
$\dim(\YY)\le 2(2k-\nu)+2d+n-\mu-5$.
If the general element of $\YY$ corresponds to a smooth $d$-gonal 
curve $X$ with a base-point-free $g^1_{d+n-\mu}$ which is not composed
with a rational involution, this inequality follows directly from what 
we have said above. Let us suppose then that a general point of
$\YY$ corresponds to a smooth $d$-gonal curve together with 
base-point-free $g^1_{d+n-\mu}$'s which is composed
with a rational involution of given degree $\gamma\ge 2$. 
Then there exists $0\le m<n$ such that $d+n-\mu=\gamma (d+m)$, and 
there exists also a dominant map 
$\ZZ\times G^1_\gamma(\mathbb{P}^1)\rightarrow\YY$
where $\ZZ$ is an irreducible component of $\GG^1_d\times_S\GG^1_{d+m}$
whose general member is a $d$-gonal curve endowed with
a base-point-free $g^1_{d+m}$ which is not composed with a rational
involution. This morphism maps a $4$-tuple $(X,f,f^\prime,\lambda)$
with $f$ a $g^1_d$ on $X$, $f^\prime$ a $g^1_{d+m}$ on $X$, and
$\lambda$ a covering $\mathbb{P}^1\rightarrow\mathbb{P}^1$ of degree $\gamma$,
to $(X,f,\lambda\circ f^\prime)$.
Again, $\ZZ$ dominates $\GG^1_d$, which implies
$\dim(\ZZ)\le 2(2k-\nu)+2d+m-5$, and furthermore
$\dim(\YY)\le 2(2k-\nu)+2d+m-5+2\gamma-2\le 2(2k-\nu)+2d+n-\mu-5$.

\medskip

We consider next $(\nu+1)$ pairs of distinct points $(x_i,y_i)$, with 
$0\le i\le \nu$ such that for any choice of $(n+1)$ pairs among them, 
$(x_{i_j},y_{i_j})$, with $0\le j\le n$ and $0\le n\le\nu$,
there exists no $L_n\in W^1_{d+n}(X)$ such that 
$h^0(X,L_n(-x_{i_j}-y_{i_j}))\ge 1$ for all $0\le j\le n$.
The $(\nu+1)$-tuple of cycles $(x_0+y_0,\dots,x_\nu+y_\nu)$
can be chosen to be generic in the space 
$X^{(2)}\times\dots\times X^{(2)}$. This is allowed by
the easy observation that, for any $n$, the incidence variety
$$\{(x_0+y_0,\dots,x_n+y_n,L_n),
h^0(L_n(-x_i-y_i))\ge 1\mbox{ for all } i\}$$ 
is at most $(2n+1)$-dimensional, whereas
$\dim\left(X^{(2)}\times\dots\times X^{(2)}\right)=2n+2$. Indeed, 
the incidence variety in question is covered by the similar incidence 
variety $\Xi$ inside $X^{2n+2}\times W^1_{d+n}(X)$,
and the fibers of the projection map from $\Xi$
to $X^{n+1}$ obtained by erasing the bundle $L_n$ and the
$y$'s are finite covers of $W^1_{d+n}(X)$ via the restriction 
of the canonical projection.

\medskip

Then we construct an irreducible curve $Y$ obtained by 
identifying $x_i$ to $y_i$ for all $i$, and denote by 
$p_i$ the corresponding node of $Y$, and by  $f:X\rightarrow Y$ the 
normalization morphism. From the genericity of the cycles $x_i+y_i$, 
the curve $Y$ can be considered to be free from non-trivial automorphisms 
and with very ample canonical bundle, {\em apply} \cite[Theorem F]{C}, 
\cite[Theorem 3.6]{CFHR}. 

We prove first that $K_{k,1}(Y,\omega_Y)=0$. 
Suppose that $K_{k,1}(Y,\omega_Y)\ne 0$.  From Proposition \ref{prop: F},
we obtain a torsion-free sheaf $F$ of rank one on $Y$
with $\chi(F)=1-k$, and $h^0(F)\ge 2$. The sheaf $F$ is either a 
line bundle, or the direct image of a line bundle on a partial 
normalization of $Y$. Observe that this partial normalization
cannot be $X$ itself. Indeed, if $F=f_*L$ with $L$ a line bundle
on $X$, then $\chi(L)=\chi(F)=1-k$, and $h^0(L)=h^0(F)\ge 2$,
which means that $L$ is a $g^1_{d-1}$ on $X$, contradicting
the hypothesis. Let us consider 
then $\varphi:Z\rightarrow Y$ the normalization of the $(\nu-n)$ 
points $p_{n+1},\dots,p_\nu$, for some $0\le n\le \nu$.
Let furthermore $\psi:X\rightarrow Z$ be the normalization of the 
remaining $(n+1)$ points $p_0,\dots,p_n$, and suppose 
$F=\varphi_*L$, for a line bundle $L$ on $Z$. 
Under these assumptions, we obtain $\chi(L)=\chi(F)=1-k$, and
so $\chi(\psi^*L)=2-k+n$, which implies that $\deg(\psi^*L)=d+n$. 
Beside, $\psi^*L$ has at least two independent sections.
Since for any node $p_i$ with $0\le i\le n$ there is a non-zero 
section of $F$ vanishing at $p_i$, it follows that 
$h^0(X,(\psi^*L)(-x_i-y_i))\ge 1$ for all $0\le i\le n$,
which contradicts the choice we made.

\medskip

We proved $K_{k,1}(Y,\omega_Y)=0$. To conclude, we apply 
\cite[Lemma 2.3]{AV} and Voisin's Remark \cite[p. 367]{V2}, 
and have $K_{k,1}(X,K_X)\subset K_{k,1}(X,K_X+x_i+y_i)
\subset K_{k,1}(Y,\omega_Y)$, for all $i$. We obtain thence 
the vanishing of $K_{k,1}(X,K_X)$ and of $K_{k,1}(X,K_X+x_i+y_i)$, 
for all $i$. The vanishing $K_{k,1}(X,K_X)=0$ is the statement of 
the Green conjecture for $X$, the fact that $\Cliff(X)$ equals 
$d-2$ being implied by the Green-Lazarsfeld non-vanishing theorem
{\em cf.} \cite[Appendix]{G}. The vanishing $K_{k,1}(X,K_X+x_i+y_i)=0$ is 
precisely the one predicted by the Green-Lazarsfeld conjecture 
for the bundle $K_X+x_i+y_i$. Then apply Theorem \ref{thm: Ap} to
conclude that the Green-Lazarsfeld conjecture is 
verified for any line bundle of degree at least $3g_X$ on $X$, 
{\em compare to} \cite[Remark 1]{AV}.
\endproof

For small $d$, one can employ classical results due H. Martens, 
Mumford and Keem on the dimensions of the Brill-Noether loci, {\em cf.} 
\cite{ACGH}, \cite{HMa}, \cite{Mu}, \cite{Ke} to obtain the following.

\begin{cor}
\label{cor: left out}
Let $X$ be a non-hyperelliptic smooth curve of gonality $d\le 6$, with 
$d<[g_X/2]+2$, and suppose that $X$ is not one of the following: 
plane curve, bielliptic, triple cover of an elliptic curve, double
cover of a curve of genus three, hexagonal curve of genus $10$ or $11$. 
Then $\Cliff(X)=d-2$ and $X$ verifies both Green, and Green-Lazarsfeld 
conjectures.
\end{cor}

\proof For a trigonal curve $X$, one has to prove that 
$\dim(W^1_{n+3}(X))\le n$ for all $0\le n\le g_X-4$.
This follows from \cite[Theorem 1]{HMa}, as we know that 
$\dim(W^1_{n+3}(X))\le n+1$ and the equality is never achieved, 
since $X$ is non-hyperelliptic.
If $d=4$, one has to prove $\dim(W^1_{n+4}(X))\le n$ 
for all $0\le n\le g_X-6$. In this case, we apply Mumford's
refinement to the Theorem of H. Martens, {\em cf.} \cite{Mu},
which shows that $\dim(W^1_{n+4}(X))\le n+1$, and equality could 
eventually hold only for trigonal (which we excluded), bielliptic 
curves or smooth plane quintics. The other cases $d=5$ and
$d=6$ are similar, and follow from \cite[Theorem 2.1]{Ke}, 
and \cite[Theorem 3.1]{Ke}, respectively.
\endproof

Some cases in the Corollary \ref{cor: left out} were known before, 
others are new. For trigonal curves, Green's conjecture was known to 
hold from Enriques and Petri, and the Green-Lazarsfeld conjecture was
verified by Ehbauer \cite{Eh}. For tetragonal curves, we knew that the 
Green conjecture was valid, {\em cf.} \cite{S2}, and \cite{V1}. All the 
other cases seem to be new. Plane curves, which were excepted from our 
statement, also verify the two conjectures, {\em cf.} \cite{Lo}, and 
\cite{A1}. Note that in a number of other cases for which our result 
does not apply, Green's conjecture is nonetheless satisfied, for
instance, for hexagonal curves of genus $10$ and Clifford index 
$3$, complete intersections of two cubics in $\mathbb{P}^3$, {\em see} 
\cite{Lo}.

\medskip

For large $d$ we cannot give similar precise results, but we still
obtain a number of examples for which Theorem \ref{thm: main1}
can be applied. For instance, curves of even genus which are 
Brill-Noether-Petri generic satisfy the hypothesis of Theorem 
\ref{thm: main1}, so they verify the two syzygy conjectures.
Other cases are obtained by looking at curves on some surfaces, 
when the special geometry of the pair $(curve,surface)$ is used,
as in the following.

\begin{cor}
\label{cor: k3} 
Let $X$ be a smooth curve of genus $2k$ and maximal Clifford index $k-1$, 
with $k\ge 2$ abstractly embedded in a $K3$ surface. Then $X$ verifies 
the Green conjecture.
\end{cor}

\proof
Since the Clifford index of $X$ is maximal, and Clifford index is 
constant in the linear system of $X$, {\em cf.} \cite{GL2}, the gonality is 
also maximal, and thus constant for smooth curves in $|X|$. Then  
the hypotheses of \cite[Lemma 3.2 (b)]{CP} are verified, which implies 
that a general smooth curve in the linear system $|X|$ has only finitely
many pencils of degree $k+1$. From Theorem \ref{thm: main1} it follows 
that the Green conjecture is verified for a general smooth curve 
$C\in |X|$, that is $K_{k,1}(C,K_C)=0$. By applying Green's hyperplane 
section theorem \cite[Theorem (3.b.7)]{G} twice, we obtain 
$K_{k,1}(X,K_X)=0$, which means that $X$ satisfies Green's conjecture, 
too. 
\endproof

Note that
Corollary \ref{cor: k3} does not apply to the particular curves 
considered by Voisin in \cite{V2}, \cite{V3}, as they are 
implicitly used in the proof.

\medskip

In view of Theorem \ref{thm: main1} and \cite[Statement (T)]{HMa}, it 
seems that understanding the geometry of curves which carry infinitely 
many minimal pencils plays a crucial role in the quest for a complete 
solution to the two conjectures. The problem of studying these 
curves has already been raised by Eisenbud, Lange, G. Martens and 
Schreyer, \cite[Remark 3.8]{ELMS}. For the beginning, it would be 
interesting to know whether the locus of smooth curves of even genus $2k$ 
with maximal Clifford index and infinitely many $g^1_{k+1}$'s is non-empty.

\section{Proof of Theorem \ref{thm: main2}}

\proofodd
We start by noting that $\dim(W^1_{k+1}(X))=1$, {\em see} \cite{FHL},  
\cite[Lemma IV.(3.3) p. 181 and Ex. VII.C-2, p. 329]{ACGH}. Then one can 
find three distinct points $x$, $y$ and $z$ of $X$ which do not belong at 
the same time to a pencil of degree $k+1$. As in the previous proof,
the cycle $x+y+z$ can be generically chosen in $X^{(3)}$, since
in our case the incidence variety 
$\{(x+y+z,L)\in X^{(3)}\times W^1_{k+1}(X),h^0(L(-x-y-z))\ge 1\}$
is two-dimensional, so that the image of its projection to $X^{(3)}$
is a surface.

For these three points, we shall prove that $K_{k,1}(X,K_X+x+y+z)=0$, 
which will imply, by means of Theorem \ref{thm: Ap}, that the 
Green-Lazarsfeld conjecture is verified for any line bundle of degree at 
least $3g_X+1$ on $X$, {\em compare to} \cite[Remark 2.6]{A2}.

\medskip

We suppose to the contrary that $K_{k,1}(X,K_X+x+y+z)\ne 0$, and shall 
reach a contradiction. For this aim, we adapt the arguments already used 
in the proof of Theorem \ref{thm: main1} to this new situation. We 
introduce a curve $Y$ with two irreducible components: the first one is 
$X$, and the second one is a smooth rational curve $E$ which passes 
through the points $x$, $y$ and $z$. Then $Y$ is a stable curve of 
arithmetic genus $g=2k+1$, and $K_{k,1}(Y,\omega_Y)\cong 
K_{k,1}(X,K_X+x+y+z)$, {\em compare to}
\cite[Proof of Lemma 2.4]{A2}. As before, from the genericity of the cycle 
$x+y+z$, we can suppose $Y$ free from non-trivial automorphisms and 
with very ample canonical bundle, {\em cf.} \cite[Theorem F]{C},
\cite[Theorem 3.6]{CFHR}.

\medskip

From Proposition \ref{prop: F}, we obtain a torsion-free,
$\omega_Y$-semistable sheaf $F$ of rank one on $Y$ with $\chi(F)=1-k$ 
and $h^0(Y,F)\ge 2$. 
We show that $F$ yields either to a pencil of degree $k+1$ on $X$ 
which passes through $x$, $y$, and $z$ or to a pencil of degree at 
most $k$. Let $F_E$, and $F_X$ be the torsion-free parts
of the restrictions of $F$ to $E$ and $X$, respectively.
It is well-known that there is a natural injection
$F\rightarrow F_E\oplus F_X$ whose cokernel is supported at the 
points among $x$, $y$ and $z$ where $F$ is invertible.
We distinguish now four cases according to the
number of nodes where $F$ is invertible.

\begin{enumerate}
\item[(i)]
$F$ is invertible at all three $x$, $y$, and $z$. 
\end{enumerate}
In this case,
$F_{|E}=F_E$, $F_{|X}=F_X$, and we have two exact sequences
$$
0\rightarrow F_E(-3)
\rightarrow F\rightarrow F_X\rightarrow 0,
$$
and, respectively,
$$
0\rightarrow F_X(-x-y-z)
\rightarrow F\rightarrow F_E\rightarrow 0.
$$

The subsheaves $F_E(-3)$ and $F_X(-x-y-z)$
are of multiranks $(1,0)$ and, respectively $(0,1)$, and,
since $\deg(\omega_{Y|E})=1$, and 
$\deg(\omega_{Y|X})=2g_X+1$, 
their  $\omega_Y$-slopes are equal to
$$\mu(F_E(-3))=\chi(F_E(-3))=\deg(F_E)-2,$$ 
and, respectively, 
$$\mu(F_X(-x-y-z))=\frac{\chi(F_X(-x-y-z))}{2g_X+1}
=\frac{\deg(F_X)-2-g_X}{2g_X+1},$$
see, for example \cite[Definition 1.1]{P}.
The $\omega_Y$-slope of $F$ equals
$$\mu(F)=\frac{1-k}{2g_X+2}.$$
From the $\omega_Y$-semistability of $F$, we
obtain $\deg(F_E)-2\le \mu(F)<0,$
which implies $\deg(F_E)\le 1$, and also
$(\deg(F_X)-1-2k)/(4k-1)\le (1-k)/(4k)$, 
which shows that $\deg(F_X)\le k+2$.

The first exact sequence implies
$\chi(F_X)=\chi(F)-\chi(F_E(-3))
=3-k-\deg(F_E)$, so $\deg(F_X)=k+1-\deg(F_E)$.

Let us suppose $\deg(F_X)=k+2$, which implies $F_E\cong \OO_E(-1)$.
Then any global section of $F$ vanishes along $E$, and
thus it vanishes at all the three points $x$, $y$, and $z$.
Since $F$ has at least two sections, the sublinear
system $H^0(F)\subset H^0(F_X)$ on $X$ has $x$,
$y$, and $z$ as base-points, in particular $h^0(F_X(-x-y-z))\ge 2$.
Then $X$ would carry a $g^1_{k-1}$, fact which contradicts 
the hypothesis.

Supposing now $\deg(F_X)\le k+1$, from the first exact sequence 
we obtain $h^0(X,F_X)\ge h^0(Y,F)\ge 2$, and since 
$X$ does not carry a $g^1_k$, it follows that $F_E=\OO_E$, and
$F_X$ is a base-point-free $g^1_{k+1}$ on $X$.
Let $\sigma$ be a non-zero global section of $F_X$ which vanishes at $x$;
such a $\sigma$ exists as $h^0(X,F_X)=2$. Then $\sigma$
is the restriction of global section $\sigma_0$ of $F$,
as the restriction morphism on global sections is 
in this case an isomorphism. 
Since $F_E=\OO_E$, the restriction of $\sigma_0$ to $E$ is a 
constant function. But $\sigma_0$ vanishes at $x$, so it vanishes 
on the whole $E$. In particular, $\sigma_0$ vanishes at $y$ and $z$, 
and thus $\sigma$ vanishes at $y$ and $z$ as well. This is a contradiction,
as we supposed that there was no such a $\sigma$.

\medskip

\begin{enumerate}
\item[(ii)]
$F$ is invertible at $y$, and $z$, and is not
invertible at $x$. 
\end{enumerate}
Let $f:Z\rightarrow Y$ be the normalization
of the point $x$. Then $F$ is the direct image of a 
line bundle $L$ on $Z$. Observe that
$\chi(Z,L) = \chi(Y,F)=1-k$,
and $h^0(Z,L) = h^0(Y,F)\ge 2$. 
We consider the exact sequences
$$
0\rightarrow F_E(-2)\rightarrow L\rightarrow F_X\rightarrow 0,
$$
and, respectively,
$$
0\rightarrow F_X(-y-z)\rightarrow L\rightarrow F_E\rightarrow 0,
$$
push them forward on $Y$, and argue as before. Since
$\chi(F_E(-2))=\chi(f_*(F_E(-2)))$
(and the same for the other sheaf), we deduce, via the
semistability of $F$, that
$\deg(F_E)\le 0$, and $\deg(F_X)\le k+1$.
The relation between the degrees of $F_E$ and $F_X$
is now $\deg(F_X)=k-\deg(F_E)$. Since $F_E$
can not be trivial (this would imply that
$F_X$ is a $g^1_k$), this shows that
$F_E$ can only be equal to $\OO_E(-1)$, 
and $F_X$ is a $g^1_{k+1}$. Then any 
global section of $L$ will be identically zero
on $E$, so it will vanish at $y$ and $z$. By the assumption on
the gonality we obtain $h^0(Z,L)=h^0(X,F_X)=2$,
which implies that any non-zero global section
of $F_X$ vanishes at $y$, and $z$, so 
$y$ and $z$ are base points of the linear system 
$|F_X|$. Then $F_X(-y-z)$ is a $g^1_{k-1}$ on $X$,
which contradicts the hypothesis. 

\begin{enumerate}
\item[(iii)]
$F$ is invertible at $x$, and is not
invertible at $y$ and $z$. 
\end{enumerate}
Is similar to the previous case.
Let $f:Z\rightarrow Y$ be the normalization
of the points $y$ and $z$. Then $F$ is the direct image of a 
line bundle $L$ on $Z$ with
$\chi(L) = \chi(F)=1-k$,
and $h^0(Z,L) = h^0(Y,F)\ge 2$. 
We consider the exact sequences
$$
0\rightarrow F_E(-1)\rightarrow L\rightarrow F_X\rightarrow 0,
$$
and, respectively,
$$
0\rightarrow F_X(-x)\rightarrow L\rightarrow F_E\rightarrow 0,
$$
and prove that $\deg(F_E)\le -1$, and $\deg(F_X)\le k$.
Then $F_X$ is a $g^1_k$ on $X$, contradiction.

\begin{enumerate}
\item[(iv)]
$F\cong F_E\oplus F_X$. 
\end{enumerate}
In this case, since $F$ is semistable, 
on the one hand, we have $\chi(F_E)\le \chi(F)/(4k)$, 
and $\chi(F_X)\le \chi(F)(4k-1)/(4k)$, and on the other hand
$\chi(F_E)+\chi(F_X)=\chi(F)$. This implies that
the two inequalities above must be equalities, which
is actually impossible, since $(1-k)/(4k)$ is not an integer
if $k\ge 2$.
\endproof

\begin{rmk}
{\rm 
Stable curves with one smooth rational irreducible component
have already been used in connection with the Green conjecture,
{\em cf.} \cite[Proof of Theorem 4]{Ei}, but with a somewhat different
purpose.
}
\end{rmk}

\section{Acknowledgements}

The present work emerged from some 
enlightening discussions with Claire Voisin, to whom I address my warmest
thanks. I thank Fabrizio Catanese for having made pertinent remarks on the 
manuscript, and indicating me the very-ampleness criteria for the
canonical bundle of Gorenstein curves. Thanks also to Lucia Caporaso, 
Gavril Farkas, Frank Schreyer, Monserrat Teixidor i Bigas, for 
useful discutions on the topic, to the IHES Bures-sur-Yvette, where this 
work was started, and the University of Bayreuth, especially Thomas 
Peternell, for hospitality, and the Alexander von Humboldt Foundation for 
financial support.

\medskip

\small

{\sc Institute of Mathematics "Simion Stoilow"
of the Romanian Academy
P.O. Box 1-764 RO-014700 Bucharest Romania}

{\em E-mail address:}
{\tt Marian.Aprodu@imar.ro}

\medskip

{\sc Universit\"at Bayreuth, Mathematisches 
Institut, Lehrstuhl Mathematik 1, D-95447 Bayreuth, Germany}

{\em E-mail address:}
{\tt aprodu@btm8x5.mat.uni-bayreuth.de}

\end{document}